\documentclass[12pt]{article}

\usepackage{latexsym,amssymb,amsmath}

\font\defi=cmr7
\newcommand{\adef}{\,\,\textrm{\raisebox{1.5pt}{\defi :}=}\,}

\newcommand{\seno}{\textrm{sin}}

\newtheorem{lemma}{Lemma}

\newcommand{\zz}{\mathbb Z}
\newcommand{\nn}{\mathbb N}
\newcommand{\cc}{\mathbb C}
\newcommand{\modulo}[1]{\left\vert{#1}\right\vert}
\newcommand{\norma}[1]{\left\Vert{#1}\right\Vert}
\newcommand{\pent}[1]{\left[{#1}\right]}
\newcommand{\pfrac}[1]{{#1}-\pent{{#1}}}

\begin{document}

\begin{center}

{\bf\Large Fourier Coefficients of Beurling Functions and a Class of
Mellin Transform Formally Determined by its Values on the Even
Integers}

\vspace{5mm}

F. Auil

\end{center}

\vspace{0.2cm}


\begin{quote}{\small
{\bf Abstract.} It is a well-known fact that 
Riemann Hypothesis will follows if the function identically equal to 
$-1$ can be arbitrarily approximated in the norm $\norma{.}$ of $L^{2}([0,1],dx)$
by functions of the form 
$f(x)=\sum_{k=1}^{n}a_{k}\,\rho\left(\frac{\theta_{k}}{x}\right)$, 
where $\rho(x)\adef\pfrac{x}$, and $a_{k}\in\cc$, 
$0<\theta_{k}\leqslant 1$ satisfies 
$\sum_{k=1}^{n}a_{k}\,\theta_{k}=0$. 
Parsevall Identity 
$\norma{f(x)+1}^{2}=\sum_{n\in\zz}\modulo{c(n)}^{2}$ is a possible
tool to compute or estimate this norm. In this note we give
an expression for the Fourier coefficients $c(n)$ of $f+1$, when $f$ is a
function defined as above. As an application, we derive an expression for 
$M_{f}(s)\adef\int_{0}^{1}(f(x)+1)\,x^{s-1}\,dx$ as a series that only
depends on $M_{f}(2k)$, $k\in\nn$. We remark that 
the Fourier coefficients $c(n)$ depend on $M_{f}(2k)$ which, for a
function $f$ defined as above, can be
expressed also in terms of the $a_{k}$'s and
$\theta_{k}$'s. Therefore, a better control on these parameters will
allow to estimate $M_{f}(2k)$ and therefore eventually to handle 
$\norma{f+1}$ via our expression for the Fourier coefficients
and Parsevall Identity. 
}
\end{quote}


\section{Introduction}

Denote as $\pent{x}$ the ``integer part of $x$'', i.e. the greatest
integer less than or equal to $x$ and define the ``fractional part''
function by $\rho(x)=\pfrac{x}$. When $f$ is a function of the form
\begin{equation}
f_{N}(x)\adef\sum_{k=1}^{N} a_{k}\,\rho\left(\frac{\theta_{k}}{x}\right), 
\label{eq-1}
\end{equation} 

\noindent
where $N\in\nn$, $a_{k}\in\cc$ and $0<\theta_{k}\leqslant 1$ for all $k\in\nn$,
an elementary computation shows 
\begin{equation}
\int_{0}^{1} (f_{N}(x)+1) \, x^{s-1} \, dx = \frac{\sum_{k=1}^{N} a_{k}\,\theta_{k}}{s-1}
+ \frac{1}{s}\left(1-\zeta(s) \, \sum_{k=1}^{N}
  a_{k}\,\theta_{k}^{s}\right);
\label{eq-2}
\end{equation}

\noindent
for $\sigma>0$, where, as usual, we denote $s=\sigma+it$. See, for
instance, \cite[p. 253]{don} for a proof. 
We will assume that function $f_{N}$ (\ref{eq-1}) satisfies also
the additional condition 
\begin{equation}
\sum_{k=1}^{N} a_{k}\,\theta_{k}=0. 
\label{eq-3}
\end{equation}

Identity (\ref{eq-2}) is the starting point of a theorem by
Beurling. See \cite[p. 252]{don} for a proof and further references.
As in the proof of (the easy half of) Beurling's Theorem, application of Schwarz
inequality to the left side of (\ref{eq-2}) allows to show that a
sufficient condition for Riemann Hypothesis is that $\norma{f(x)+1}$ be
done arbitrarily small for a convenient choice of $a_{k}$ and
$\theta_{k}$, where $\norma{.}$ denotes the norm in
$L^{2}([0,1],dx)$.\\

Just for reference, we will call a function $f_{N}$ as in (\ref{eq-1})
as an \emph{approximation} or \emph{Beurling function}, and the sequence
$\{f_{N}\}_{N\in\nn}$ is called an \emph{approximation
sequence}. We remark that approximation sequences do not necessarily
\emph{converge} to $-1$ in $L^{2}([0,1],dx)$; see the excellent work
\cite{bae} on this topic.
For abuse of notation and language, when $f_{N}$ is a Beurling
function we will call $f_{N}+1$ also a Beurling function.\\

A method to compute 
$\norma{f_{N}+1}\adef\left(\int_{0}^{1}\modulo{f_{N}(t)+1}^{2}\,dt\right)^{1/2}$
would be to use not this definition but the Parsevall Identity 
$\norma{f_{N}+1}^{2}=\sum_{n\in\zz}\modulo{c(N,n)}^{2}$. In the
Sec. \ref{fourier} we give an expression for the Fourier coefficients
$c(N,n)$ of the Beurling function $f_{N}+1$.\\

In this note, unless explicit statement on
contrary, we assume condition (\ref{eq-3}) on a Beurling function. 
At a certain point we will assume also that $\theta_{k}=1/b_{k}$ 
where $b_{k}\in\nn$ and $\modulo{a_{k}}\leqslant 1$. This restriction
on the $\theta_{k}$'s is not serious at all after Theorem 1.1 in
\cite{bae2}. On the other hand, the restriction on the $a_{k}$'s
include some of the so-called \emph{natural approximations} considered
in \cite{bae}.


\section{The Fourier Coefficients for a Beurling function}
\label{fourier}

For convenience we define $F_{N}\adef f_{N}(x)+1$. 
We extend $F_{N}$ to an \emph{odd} function in $[-1,1]$. Therefore, 
$\int_{-1}^{1}F_{N}(x)\,\cos(n\pi x)\,dx = 0$, and

\begin{multline*}
c(N,n)\adef\int_{-1}^{1}F_{N}(x)\,\seno(n\pi x)\,dx 
= 2\,\int_{0}^{1}F_{N}(x)\,\seno(n\pi x)\,dx\\
= 2\left[\int_{0}^{1}\seno(n\pi x)\,dx 
+ \sum_{k=1}^{N}a_{k}\,\int_{0}^{1}\rho\left(\frac{\theta_{k}}{x}\right)\,\seno(n\pi
x)\,dx \right]\\
= 2\left[\left.-\frac{\cos(n\pi x)}{n\pi}\right\vert_{0}^{1}
+
\sum_{k=1}^{N}a_{k}\left(\sum_{j=1}^{\infty}\int_{\frac{\theta_{k}}{j+1}}^{\frac{\theta_{k}}{j}}\rho\left(\frac{\theta_{k}}{x}\right)\,\seno(n\pi
x)\,dx 
+ \int_{\theta_{k}}^{1}\rho\left(\frac{\theta_{k}}{x}\right)\,\seno(n\pi
x)\,dx \right)\right]
\end{multline*}
\begin{multline*}
= 2\left[\frac{1}{n\pi}(1-\cos(n\pi))\right.\\
+
\left.\sum_{k=1}^{N}a_{k}\left(\sum_{j=1}^{\infty}\theta_{k}\,\int_{\frac{\theta_{k}}{j+1}}^{\frac{\theta_{k}}{j}}\frac{\seno(n\pi
x)}{x}\,dx - j\,\int_{\frac{\theta_{k}}{j+1}}^{\frac{\theta_{k}}{j}}\seno(n\pi
x)\,dx 
+ \theta_{k}\,\int_{\theta_{k}}^{1}\frac{\seno(n\pi x)}{x}\,dx \right)\right]\\
= \frac{2}{n\pi}(1-\cos(n\pi))\;
+\;
2\,\sum_{k=1}^{N}a_{k}\left(\theta_{k}\,\int_{0}^{1}\frac{\seno(n\pi x)}{x}\,dx
- 
\sum_{j=1}^{\infty}j\,\int_{\frac{\theta_{k}}{j+1}}^{\frac{\theta_{k}}{j}}\seno(n\pi
x)\,dx \right)\\
= \frac{2}{n\pi}(1-\cos(n\pi))
+ 2\left(\sum_{k=1}^{N}a_{k}\,\theta_{k}\right)\left(\int_{0}^{1}\frac{\seno(n\pi x)}{x}\,dx\right)
- 2\,\sum_{k=1}^{N}a_{k}\,\sum_{j=1}^{\infty}j\,\left.\frac{-\cos(n\pi
x)}{n\pi}\right\vert_{\frac{\theta_{k}}{j+1}}^{\frac{\theta_{k}}{j}}
\end{multline*}

Denote the first term in last line by
$A_{1}\adef\frac{2}{n\pi}(1-\cos(n\pi))$, and note that the second
term vanishes because (\ref{eq-3}), hence

\begin{equation}
c(N,n) = A_{1}
+ \frac{2}{n\pi}\,\sum_{k=1}^{N}a_{k}\,\sum_{j=1}^{\infty}j\,\left[\cos\left(\frac{n\pi\theta_{k}}{j}\right)
- \cos\left(\frac{n\pi\theta_{k}}{j+1}\right)\right].
\label{eq10}
\end{equation}

Now, if $L$ is any natural number, replacing in (\ref{eq10}) the expression for the
$L$-th Taylor approximation for $\cos x$ given by 

\begin{equation}
\cos x = \sum_{l=0}^{L} (-1)^{l}\,\frac{x^{2l}}{(2l)!} \; 
+ \;\frac{1}{L!}\,\int_{0}^{x}\cos^{(L+1)}(t)\,(x-t)^{L}\,dt,
\end{equation}

\noindent
we get
\begin{multline}
c(N,n)\\
= A_{1}
+ \frac{2}{n\pi}\,\sum_{k=1}^{N}a_{k}\,\sum_{j=1}^{\infty}j\,\left(\sum_{l=0}^{L}(-1)^{l}\,\frac{(n\pi\theta_{k})^{2l}}{(2l)!}
\left[\frac{1}{j^{2l}}-\frac{1}{(j+1)^{2l}}\right] \; + \;
R(L,j,n,k)\right)\\
= A_{1}
+ \frac{2}{n\pi}\,\sum_{k=1}^{N}a_{k}\left(\sum_{l=1}^{L}(-1)^{l}\,\frac{(n\pi\theta_{k})^{2l}}{(2l)!}\,
\sum_{j=1}^{\infty}j\,\left[\frac{1}{j^{2l}}-\frac{1}{(j+1)^{2l}}\right]
\; 
+ \; \sum_{j=1}^{\infty}j\,R(L,j,n,k)\right);
\label{eq-20}
\end{multline}
where
\begin{multline}
R(L,j,n,k)\adef
\frac{1}{L!}\,\int_{0}^{\frac{n\pi\theta_{k}}{j}}\cos^{(L+1)}(t)\,\left(\frac{n\pi\theta_{k}}{j}-t\right)^{L}\,dt\\
-
\;\frac{1}{L!}\,\int_{0}^{\frac{n\pi\theta_{k}}{j+1}}\cos^{(L+1)}(t)\,\left(\frac{n\pi\theta_{k}}{j+1}-t\right)^{L}\,dt.
\label{eq-21}
\end{multline}
Now, observe that
\begin{multline}
\sum_{j=1}^{\infty}j\,\left[\frac{1}{j^{2l}}-\frac{1}{(j+1)^{2l}}\right]
=\lim_{J\to\infty}\sum_{j=1}^{J}\left[\frac{j}{j^{2l}}-\frac{j+1-1}{(j+1)^{2l}}\right]\\
=\lim_{J\to\infty}\left(\sum_{j=1}^{J}\left[\frac{j}{j^{2l}}-\frac{j+1}{(j+1)^{2l}}\right]
\;+\;\sum_{j=1}^{J}\frac{1}{(j+1)^{2l}}\right)\\
=\lim_{J\to\infty}\left(1\; - \;\frac{J+1}{(J+1)^{2l}} \; + \; \sum_{j=2}^{J+1}\frac{1}{j^{2l}}\right)
= 1\;+\;\sum_{j=2}^{\infty}\frac{1}{j^{2l}} = \zeta(2l);
\label{eq-30}
\end{multline}
Substituting (\ref{eq-30}) in (\ref{eq-20}) and 
denoting $R(L,n,k)\adef \sum_{j=1}^{\infty}j\,R(L,j,n,k)$, we get
\begin{multline}
c(N,n)
= A_{1}
+ \frac{2}{n\pi}\,\sum_{k=1}^{N}a_{k}\left(\sum_{l=1}^{L}(-1)^{l}\,\frac{(n\pi\theta_{k})^{2l}}{(2l)!}\,
\zeta(2l)\; 
+ \; R(L,n,k)\right)\\
= A_{1}
+
\frac{2}{n\pi}\,\sum_{l=1}^{L}(-1)^{l}\,\frac{(n\pi)^{2l}}{(2l)!}\,\zeta(2l)
\left(\sum_{k=1}^{N}a_{k}\,\theta_{k}^{2l}\right)\;
+ \; \sum_{k=1}^{N}a_{k}\,R(L,n,k).
\label{eq-40}
\end{multline}

If we denote $M_{g}(s)\adef\int_{0}^{1} g(x)\,x^{s-1}\,dx$, then by
(\ref{eq-2}), under condition (\ref{eq-3}), we have 

\begin{equation}
\sum_{k=1}^{N}a_{k}\,\theta_{k}^{2l}
= \frac{1}{\zeta(2l)}\left(1 \;-\;(2l)\,M_{F_{N}}(2l)\right).
\label{eq-50}
\end{equation}

Substituting now (\ref{eq-50}) in (\ref{eq-40}) we get

\begin{multline}
c(N,n)
= A_{1}\;
+ \; \frac{2}{n\pi}\,\sum_{l=1}^{L}(-1)^{l}\,\frac{(n\pi)^{2l}}{(2l)!}\\
- \; \frac{2}{n\pi}\,\sum_{l=1}^{L}(-1)^{l}\,\frac{(n\pi)^{2l}}{(2l)!}\,(2l)\,M_{F_{N}}(2l)
+ \; \sum_{k=1}^{N}a_{k}\,R(L,n,k)\\
\hspace{-5cm}= A_{1}\;
+ \; \frac{2}{n\pi}\,\sum_{l=1}^{L}(-1)^{l}\,\frac{(n\pi)^{2l}}{(2l)!}\\
+\; 2\sum_{l=1}^{L}(-1)^{l-1}\,\frac{(n\pi)^{2l-1}}{(2l-1)!}\,M_{F_{N}}(2l)
+ \; \sum_{k=1}^{N}a_{k}\,R(L,n,k).
\label{eq-60}
\end{multline}

This is an exact expression valid for all $L\in\nn$ but a better
expression arises in the limit $L\to\infty$. In this case, the first
term cancels the second and is not difficult to prove the following 

\bigskip

\begin{lemma} For each $k\in\nn$ assume $\modulo{a_{k}}\leqslant 1$
and $\theta_{k}=1/b_{k}$ with $b_{k}\in\nn$. Then, for any fixed $N$ and $n$
in $\nn$ is $\displaystyle{\lim_{L\to\infty}\sum_{k=1}^{N}a_{k}\,R(L,n,k)=0}$. 
\label{lemma-1}
\end{lemma}

\bigskip

\noindent
(The proof will be given in Sec. \ref{proof}). Therefore, the final expression for
the Fourier coefficient is

\begin{equation}
c(N,n) = 2\sum_{l=1}^{\infty}(-1)^{l-1}\,\frac{(n\pi)^{2l-1}}{(2l-1)!}\,M_{F_{N}}(2l).
\label{eq-70}
\end{equation}


\section{An Expression for the Mellin Transform}

Directly from (\ref{eq-70}) we have

\begin{multline}
M_{F_{N}}(s)=\int_{0}^{1}F_{N}(x)\,x^{s-1}\,dx
=\int_{0}^{1}\left(\sum_{n=1}^{\infty}c(N,n)\,\seno(n\pi x)\right)\,x^{s-1}\,dx\\
=\sum_{n=1}^{\infty}\left(\int_{0}^{1}\seno(n\pi x)\,x^{s-1}\,dx\right)\, 
2\sum_{l=1}^{\infty}(-1)^{l-1}\,\frac{(n\pi)^{2l-1}}{(2l-1)!}\,M_{F_{N}}(2l).
\end{multline}

The last expression gives the value of $M_{F_{N}}(s)$ in term of its
values in the even integers $M_{F_{N}}(2l)$.\\

Note also that under the hypothesis of Lemma \ref{lemma-1} we have 
$\modulo{\sum_{k=1}^{N}a_{k}\,\theta_{k}^{2l}}\leqslant\zeta(2l)$,
and this combined with (\ref{eq-50}), or (\ref{eq-2}), gives an
estimation on the Mellin transform
\begin{equation}
\modulo{M_{F_{N}}(2l)}\leqslant\frac{1+\zeta^{2}(2l)}{2l}.
\end{equation}

As in the derivation of relation (\ref{eq-200}) in Sec. \ref{proof},
a better control on the $a_{k}$'s and $\theta_{k}$'s will allows to control
$M_{F_{N}}(2l)$ via (\ref{eq-50}), and therefore to control
$\norma{f_{N}(x)+1}$ via
(\ref{eq-70}) and Parsevall Identity.


\section{Proof of Lemma \ref{lemma-1}}
\label{proof}

For sake of brevity, denote
$I(j)\adef(1/L!)\,\int_{0}^{\frac{n\pi\theta_{k}}{j}}\cos^{(L+1)}(t)\,(\frac{n\pi\theta_{k}}{j}-t)^{L}$.
From (\ref{eq-21}) we have
\begin{multline}
R(L,n,k)\adef\sum_{j=1}^{\infty}j\,R(L,j,n,k)
=\lim_{J\to\infty}\sum_{j=1}^{J}j\,I(j)\;-\;(j+1-1)\,I(j+1)\\
=\lim_{J\to\infty}\left(\sum_{j=1}^{J}\left[j\,I(j)\;-\;(j+1)\,I(j+1)\right] \;+\;\sum_{j=1}^{J}I(j+1)\right)\\
=\lim_{J\to\infty}\left(I(1)\;-\;(J+1)\,I(J+1) \; + \;\sum_{j=1}^{J}I(j+1)\right)\\
=\lim_{J\to\infty}\left(\sum_{j=1}^{J+1}I(j)\;-\;(J+1)\,I(J+1)\right).
\end{multline}

Using the elementary estimative
\begin{equation}
\modulo{I(j)}
\leqslant
\frac{1}{L!}\,\left.\frac{-\left(\frac{n\pi\theta_{k}}{j}-t\right)^{L+1}}{L+1}\right\vert_{0}^{\frac{n\pi\theta_{k}}{j}}
=\frac{1}{j^{L+1}}\,\frac{(n\pi\theta_{k})^{L+1}}{(L+1)!},
\end{equation}
for any $L\in\nn$ we have 
\begin{multline}
\modulo{R(L,n,k)}\leqslant
\lim_{J\to\infty}\frac{(n\pi\theta_{k})^{L+1}}{(L+1)!}
\left(\sum_{j=1}^{J+1}\frac{1}{j^{L+1}}\; +
\;\frac{1}{(J+1)^{L}}\right)
=\frac{(n\pi\theta_{k})^{L+1}}{(L+1)!}\,\zeta(L+1).
\end{multline}
Therefore, for any $N\in\nn$, 
\begin{multline}
\modulo{\sum_{k=1}^{N}a_{k}\,R(L,n,k)}\leqslant
\sum_{k=1}^{N}\modulo{R(L,n,k)}\leqslant
\frac{(n\pi)^{L+1}}{(L+1)!}\,\zeta(L+1)\sum_{k=1}^{N}\theta_{k}^{L+1}\\
=\frac{(n\pi)^{L+1}}{(L+1)!}\,\zeta(L+1)\sum_{k=1}^{N}\frac{1}{b_{k}^{L+1}}
\leqslant\frac{(n\pi)^{L+1}}{(L+1)!}\,\zeta^{2}(L+1). 
\label{eq-200}
\end{multline}
Now the assertion of Lemma \ref{lemma-1} follows observing that
$\zeta(L+1)$ remains close to $1$ for large $L$ and the first factor in
the right side of (\ref{eq-200}) goes to zero because is the $(L+1)$-th
term of the (convergent) series for $\exp(n\pi)$.



\vspace{1.5cm}

\noindent
Fernando Auil\\

\noindent
Instituto de F\'isica\\
Universidade de S\~ao Paulo\\  
Caixa Postal 66318\\  
CEP 05315-970\\       
S\~ao Paulo - SP\\            
Brasil\\

\noindent                  
E-mail: {\tt auil@fma.if.usp.br}



\begin{thebibliography}{99}


\bibitem{don} Donoghue W. F. {\it  Distributions and Fourier
    Transforms}, Academic Press, New York, 1969. 

\bibitem{bae} L. Báez-Duarte, {\it Arithmetical Aspects of Beurling's
Real Variable Reformulation of the Riemann
Hypothesis}. arXiv:math.NT/0011254v1

\bibitem{bae2} L. Báez-Duarte, {\it A Strengthening of the
Nyman-Beurling Criterion for the Riemann
Hypothesis}. arXiv:math.NT/0202141v2


\end{thebibliography}
\end{document}